\documentclass[12pt]{article}
\usepackage{amssymb}
\usepackage{graphicx}
\usepackage{latexsym}
\def\part#1{\frac{\partial\phantom{q}}{\partial#1}}

\newenvironment{rmk}{\begin{trivlist}\item[]{\bf Remark:} }
{\end{trivlist}}
\newenvironment{rmks}{\begin{trivlist}\item[]{\bf Remarks:} }
{\end{trivlist}}
\newenvironment{ex}{\begin{trivlist}\item[]{\bf Example:} }
{\end{trivlist}}
\newenvironment{prf}{\begin{trivlist}\item[]{\bf Proof:} }
{\hfill $\Box$ \end{trivlist}}

\newtheorem{thm}{Theorem}

\newtheorem{prp}[thm]{Proposition}
\newtheorem{lemma}[thm]{Lemma}

\newcommand{\lie}[1]{\mathfrak{#1}}
\def\End{\mathop{\rm End}\nolimits}
\def\Prym{\mathop{\rm Prym}\nolimits}
\def\Hom{\mathop{\rm Hom}\nolimits}

\def\ker{\mathop{\rm ker}\nolimits}
\def\coker{\mathop{\rm coker}\nolimits}

\def\Pic{\mathop{\rm Pic}\nolimits}

\def\deg{\mathop{\rm deg}\nolimits}

\def\tr{\mathop{\rm tr}\nolimits}

\newcommand{\C}{\mathbf{C}}

\newcommand{\Z}{\mathbf{Z}}

\newcommand{\PP}{{\rm P}}

\begin{document}
\title{Spinors, Lagrangians and rank $2$ Higgs bundles}
\author{Nigel Hitchin\\[5pt]}

\maketitle

\section{Introduction}
The moduli space ${\mathcal M}$ of Higgs bundles on a Riemann surface has been studied for almost 30 years now but still presents some conceptual problems. There are two standard viewpoints:  on the one hand it has an open subset which is the cotangent bundle of the compact moduli space ${\mathcal N}$ of semistable bundles, and on the other it is a fibration over a vector space with generic fibre an abelian variety.  The holomorphic symplectic structure is the common feature here, the fibres of both being Lagrangian submanifolds, the second case forming an algebraically completely integrable system.  This paper concerns another construction of Lagrangians due to D.Gaiotto, which we use  to investigate the simplest case of rank $2$ Higgs bundles. The Lagrangians appear in different ways corresponding to the two viewpoints.  

The starting point is a    complex semi-simple Lie group $G^c$ and a holomorphic principal $G^c$-bundle over a compact Riemann surface $\Sigma$,  and in addition a symplectic representation of $G^c$. The symplectic condition means that there is an invariant quadratic function $\mu$ on the representation space, the moment map, with values in the dual of the Lie algebra, or using the Killing form, in the Lie algebra itself. 
Let $W$ be the  associated  vector bundle. We take a spinor field, a holomorphic section $\psi$ of $W\otimes K^{1/2}$, where $K$ is the canonical bundle, and apply the moment map to obtain $\mu(\psi)=\Phi$ a section of $\lie{g}\otimes K$, where $\lie{g}$ denotes the adjoint bundle.  The pair $(W,\Phi)$ is a Higgs bundle and general principles suggest that, when this is stable, the set of such pairs should be Lagrangian in the Higgs bundle moduli space.

Current interest in these is driven both by the fact that these are so-called BAA-branes which should have mirror BBB-branes on the moduli space for the Langlands dual group, but also (and related) is the issue of what they represent in some holomorphic Fukaya category. In particular one might consider the morphisms to the Lagrangian torus fibres.

In this paper we consider examples: the simplest case where $G^c=SL(2,\C)$.The irreducible symplectic representations are the odd symmetric powers $S^mV$. When $V$ is  stable the condition for the existence of a holomorphic section of $S^mV\otimes K^{1/2}$ is the vanishing of the Quillen determinant section of a certain holomorphic line bundle on the  moduli space of stable bundles ${\mathcal N}$. The smooth points of this determinant divisor correspond to the case of a one-dimensional space of sections $\psi$. In this case we shall see that the Lagrangian defined by $\mu(\psi)$ can be viewed as the conormal bundle of the divisor naturally embedded in the cotangent bundle. This is the picture according to the first approach.

 From the integrable system viewpoint, when $m=1$ the Lagrangian turns out to be  a particular algebraic component of the nilpotent cone, the fibre over zero in the  integrable system. For more general $m$ the situation is more interesting. When $m=3$ the Lagrangian intersects the generic fibre of the integrable system in the non-trivial points of order $3$. Conversely, we show that for any odd $m$, an $m$-period point in a generic fibre is represented by a Higgs bundle $(S^mV,\mu(\psi))$ for  a canonical choice of spinor $\psi$.

If $\Sigma$ has genus $2$, then ${\mathcal N}$ is well-known to be the projective space $\PP^3$ with the quartic Kummer surface corresponding to the semistable locus. For $m=1$ the determinant divisor is a plane. There are $16$ of these corresponding to the $16=2^{2g}$ choices for $K^{1/2}$. These are classically known as {\it tropes}. The corresponding Lagrangian in ${\mathcal M}$ is isomorphic to $\PP^3$. 

The determinant divisor for $m=3$ is a surface of degree $10$ in $\PP^3$. We describe  its intersection with the tropes, which in fact gives the intersection of the Lagrangians for $m=1$ and $3$. These intersections are reducible plane curves: in the trope corresponding to the given choice of $K^{1/2}$ it is a double conic together with a sextic, whose equation we derive. The surface intersects each of the other 15 tropes  in a sextic and a quartic. Their intersection properties are related to the classical geometry of the $16_6$ configuration of points and planes. We also briefly consider the  Lagrangian for  $m=1$ in genus $3$, following \cite{OPP}. 
\vskip .25cm
The author wishes to thank Davide Gaiotto for pointing out this construction, Oscar Garcia-Prada for discussions  and EPSRC, ICMAT and QGM for support. 

\section{A general setting}
We describe here  Gaiotto's approach. 
There is an infinite-dimensional approach to Higgs bundles which is the background to the original construction in \cite{NH} and is useful in the present context. We start with a fixed $C^{\infty}$ principal $G$-bundle for $G$ the compact real form of $G^c$ and regard the space of all holomorphic structures as an affine space ${\mathcal A}$ with group of translations $\Omega^{01}(\Sigma, \lie{g})$. Given a representation $G^c\rightarrow GL(n,\C)$, the difference of any two $\bar\partial$-operators on sections of the associated bundle $W$ is $\rho(a)\in \Omega^{01}(\Sigma, \End W)$ where $\rho:\lie{g}\rightarrow \lie{gl}(n)$ is the Lie algebra homomorphism and $a\in \Omega^{01}(\Sigma, \lie{g})$. The cotangent bundle $T^*{\mathcal A}$ is the product ${\mathcal A}\times \Omega^{10}(\Sigma, \lie {g})$ using the integration pairing between $\Omega^{01}(\Sigma, \lie{g})$ and $\Omega^{10}(\Sigma, \lie{g})$. The group of complex gauge transformations acts symplectically with moment map $(A,\Phi)\mapsto \bar\partial_A\Phi$ and the symplectic quotient is then formally the space of equivalence classes of holomorphic structures together with holomorphic sections $\Phi$ of $\lie{g}\otimes K$. The reduction to $G$ allows one to define a real moment map which leads to a hyperk\"ahler quotient and the notion of stability but we shall not need that for this formal setting.

In two dimensions a holomorphic line bundle  $K^{1/2}$ such that $(K^{1/2})^2\cong K$ is a spin structure and the $\bar\partial$-operator $\bar\partial:\Omega^{00}(\Sigma, K^{1/2})\rightarrow \Omega^{01}(\Sigma, K^{1/2})$ is the Dirac operator \cite{NH0}. Given a holomorphic vector bundle $W$ we can define  the coupled Dirac operator $\bar\partial_A:\Omega^{00}(\Sigma, W\otimes K^{1/2})\rightarrow \Omega^{01}(\Sigma, W\otimes K^{1/2})$. If $W$ has an invariant symplectic form $\langle u,v\rangle$ then there is a natural functional 
\begin{equation}
L(A,\psi)=\int_{\Sigma}\langle \bar\partial_A\psi,\psi\rangle
\label{L}
\end{equation}
on ${\mathcal A}\times \Omega^{00}(W\otimes K^{1/2})$. The skew-symmetry of  $\langle u,v\rangle$ together with Stokes' theorem shows that this is a well-defined quadratic functional in $\psi$. 

Regard this as a family of functions on ${\mathcal A}$ parametrized by $\psi\in \Omega^{00}(\Sigma, W\otimes K^{1/2})$. For a  finite-dimensional manifold $M$  the derivative of a function $L$ on $M$ defines a Lagrangian submanifold of $T^*M$, but more generally if $L$ also depends on  $u\in  N$ the derivative relative to $M$ of the critical locus  with respect to 
$N$ is a subspace of $T^*M\times N$ which (under appropriate transversality conditions) projects to a Lagrangian submanifold of $T^*M$. In the infinite-dimensional case  above, replace $M$ by ${\mathcal A}$ and $N$ by the space $\Omega^{00}(\Sigma, W\otimes K^{1/2})$. Then the critical locus is obtained if 
$$0= 2\int_{\Sigma}\langle \bar\partial_A\psi,\dot\psi\rangle$$
for all  variations $\dot\psi$. This means $\bar\partial_A\psi=0$ so $\psi$ is holomorphic. 
The graph of the derivative with respect to $A$ is given by $(A,\Phi)$ such that for all $\dot A\in \Omega^{01}(\Sigma, \lie{g})$
$$\int_{\Sigma} (\Phi,\dot A)=\int_{\Sigma}\langle \rho(\dot A)\psi,\psi\rangle=\int_{\Sigma}( \dot A,\rho^*(\psi\otimes \psi))=\int_{\Sigma}(\mu(\psi), \dot A)$$
giving $\Phi=\mu(\psi)$ (here we are identifying the symmetric tensor product $\psi\otimes \psi$ with an element of  the Lie algebra of the symplectic group.) 
Quotienting by the complex gauge group, the expectation is that we will get a Lagrangian submanifold but   we can in principle only assert that the subspace of ${\mathcal M}$ is isotropic. 

To see the isotropic condition more directly, a first order deformation of the equations $\bar\partial_A\psi=0, \Phi=\rho^*(\psi\otimes \psi)$ is given by
$$\bar\partial_A\dot\psi+\rho(\dot A)\psi=0,\qquad \dot\Phi=\rho^*(\dot\psi\otimes \psi+\psi\otimes \dot \psi)$$
and then 
$$\int_{\Sigma}(\dot\Phi_1,\dot A_2)=\int_{\Sigma}(\rho^*(\dot\psi_1\otimes \psi+\psi\otimes \dot \psi_1),\dot A_2)=2 \int_{\Sigma}\langle\dot \psi_1,\rho(\dot A_2)\psi\rangle=-2 \int_{\Sigma}\langle\dot \psi_1,\bar\partial_A\dot\psi_2\rangle.$$
But (as in (\ref{L})) this is symmetric so 
$$\int_{\Sigma}(\dot\Phi_1,\dot A_2)-\int_{\Sigma}(\dot\Phi_2,\dot A_1)=0$$
which is the isotropy condition. In the examples the submanifold will  in fact be Lagrangian.

\begin{rmk} Note that if $\psi$ is holomorphic so is $\lambda \psi$ for $\lambda\in \C$ and if $(A,\mu(\psi))$ lies on the Lagrangian submanifold so does $(A,\lambda^2\mu(\psi))$. The $\C^*$-action $(A,\Phi)\mapsto (A,\lambda \Phi)$ on the Higgs bundle moduli space defines a vector field $X$ which is therefore tangential to the Lagrangian. This means that if $\omega$ is the symplectic form on ${\mathcal M}$ then $\omega(X,Y)=0$ for all tangent vectors $Y$ to the Lagrangian. Hence the one-form $i_X\omega$ (which  restricted to the cotangent bundle $T^*{\mathcal N}$ is the tautological 1-form) vanishes  and not just its derivative $d(i_X\omega)=2\omega$. This is a stronger condition than the Lagrangian one. For example, the only fibres in the integrable system that have this property are the components of the nilpotent cone \cite{BGL}.
\end{rmk}
\section{Determinant divisors}
Quillen \cite{Q} defined the determinant line bundle for a smooth family of $\bar\partial$-operators on a $C^{\infty}$ vector bundle over a Riemann surface $\Sigma$. Each such operator $\bar\partial_A$ defines a determinant line
\begin{equation}
L_A=(\Lambda^{top}\ker\bar\partial_A)^*\otimes \Lambda^{top}\coker \bar\partial _A
\label{Q}
\end{equation}
and as $A$ varies this is a smooth line bundle even though the dimension of the kernel and cokernel may jump. If the index of $\bar\partial_A$ is zero (which is a topological property and holds for the whole connected family) then there is a canonical determinant section which vanishes whenever there exists a non-trivial solution to $\bar\partial_A\psi=0$. If the bundle $W$ is symplectic and we take $\bar\partial_A:\Omega^{00}(W\otimes K^{1/2})\rightarrow \Omega^{01}(W\otimes K^{1/2})$ then the index is zero for Serre duality gives $H^0(\Sigma, W\otimes K^{1/2})^*\cong H^1(\Sigma, W^*\otimes K^{-1/2}\otimes K)=H^1(\Sigma, W\otimes K^{1/2})$ since $W\cong W^*$ using the symplectic structure.

 On a {\it moduli space} $M$ of such operators there may not exist a universal family over $\Sigma\times M$ but the determinant line bundle is always defined if the index is zero. 
This is because there is always a local universal bundle $W_U$ over a small enough open set $U\subset M$ but $W_U$ and $W_V$ over $U\cap V$ differ by multiplying by a constant scalar. This changes both the kernel and cokernel  by the same scalar and  if the dimensions are the same leaves $L_A$ in (\ref{Q}) unchanged. 

Thus on the moduli space of stable bundles ${\mathcal N}$ each symplectic representation defines a distinguished determinant divisor for the associated $\bar\partial$-operators on $W\otimes K^{1/2}$. Since the above construction requires the existence of a holomorphic section $\psi$, if we restrict to the open set $T^*{\mathcal N}\subset {\mathcal M}$ we obtain a subspace which projects to the determinant divisor in ${\mathcal N}$. If $\mu(\psi)$ is identically zero  then the subspace is  just  the divisor lying in the zero section which is isotropic but of less than maximal dimension.

\begin{prp} Let $L$ be the determinant line bundle on ${\mathcal N}$  for spinors coupled to a symplectic representation of  $G^c$ and $Y$ the open subset of the determinant divisor where the kernel is one-dimensional. Then the moment map $\mu$ defines a map from $L^*$ restricted to $Y$ to the conormal bundle of $Y$ as a submanifold of $T^*{\mathcal N}$.
\end{prp}
\begin{rmks} 
\noindent 1. The conormal bundle of any submanifold  $Y\subset X$ embeds in $T^*X$ as the cotangent vectors which annihilate the tangent vectors of $X$ and is a Lagrangian submanifold. Conversely any conic Lagrangian submanifold (i.e. invariant by rescaling in the fibre directions) is locally a conormal bundle if the derivative of the projection to $X$ has constant rank.

\noindent 2. If $Y$ is the nondegenerate zero set of a section $s$ of a line bundle $L$ then the derivative of $s$ on $Y$ is a well-defined homomorphism $L^*\rightarrow T^*$ which embeds the total space of $L^*$ as the conormal bundle  in $T^*X$.
\end{rmks}

\begin{prf} We assume that $\dim\ker \bar\partial_A=1$. The Serre duality above is a consequence of the   fact that the operator is formally symmetric (as in (\ref{L})) which means there is a natural isomorphism between $(\ker\bar\partial_A)^*$ and  $\coker \bar\partial _A$. It follows that the determinant line is 
$$L_A=(\Lambda^{top}\ker\bar\partial_A)^{-2}$$
and hence on $Y$ we have $L^*\cong (\ker\bar\partial_A)^2$.

Therefore  if $\bar\partial_A\psi=0$ then $\psi\otimes \psi$ spans the fibre of $L^*$ at $[A]$. Since the moment map is homogeneous of degree $2$ it gives a homomorphism $\mu(\psi)$ from $L^*$ to $T^*{\mathcal N}$. 

Let $\dot A\in \Omega^{01}(\Sigma, \lie{g})$ represent a tangent vector to $Y$ and evaluate the cotangent vector $\Phi=\mu(\psi)$ on $[\dot A]$. We obtain
$$\int_{\Sigma}(\Phi,\dot A)=\int_{\Sigma}(\rho^*(\psi\otimes\psi),\dot A)=2\int_{\Sigma}\langle \psi, \rho(\dot A)\psi\rangle=-2\int_{\Sigma}\langle\psi, \bar\partial_A\dot\psi\rangle$$
since $\bar\partial_A\dot\psi+\rho(\dot A)\psi=0$ is the tangency condition for $\dot A$. But $\psi$ is holomorphic so by Stokes' theorem we get zero. The image is therefore the conormal bundle.
\end{prf}

\begin{ex} If $G^c=Sp(2m,\C)$ and we take the defining representation then $\mu(\psi)=\psi\otimes \psi$ which is never zero. The determinant divisor here is called the generalized theta-divisor and its cohomology class generates $H^2({\mathcal N},\Z)$. The Lagrangian in ${\mathcal M}$ therefore intersects $T^*{\mathcal N}$ in the conormal bundle of this divisor.
\end{ex}

\section{Lagrangians  for $SL(2,\C)$}
\subsection{Spectral curves}\label{spec}
 We turn next to the integrable system setting, where we restrict to the case $G^c=SL(2,\C)$.
The fibration picture for this group  involves the {\it spectral curve}, a double covering of $\Sigma$ \cite{NH},\cite{BNR}. What follows is a brief account of its properties.

A Higgs bundle for the basic representation on $\C^2$ is a rank $2$ holomorphic vector bundle $V$ with a symplectic form and a trace zero Higgs field $\Phi\in H^0(\Sigma,\End V\otimes K)$. The characteristic equation $\det(x-\Phi)=x^2-q=0$ defines a curve $S$ in the total space of the canonical bundle $\pi:K\rightarrow \Sigma$. If the quadratic differential $q$ has simple zeros then $S$ is smooth. Invariantly speaking $x$ can be regarded as the tautological section of $\pi^*K$ and $\det \Phi$ is a section of $K^2$ so $S$ is a divisor in the class of $\pi^*K^2$. It follows that $K_S\cong \pi^*K^2$. The involution $x\mapsto -x$ on $K$ preserves $S$ and we denote this by $\sigma$. Its fixed point set is the zero section $x=0$ and $\Sigma$ is the quotient curve.

On $S$, $\pm x$ is an eigenvalue of $\Phi$ and  $\pi^*V$ has two well-defined subbundles $L^*$ and $\sigma^*L^*$ giving the eigenspaces. The two coincide on $x=0$ and so the homomorphism $\sigma^*L^*\subset \pi^*V\rightarrow \pi^*V/L^*=L$ shows that $(\sigma^*L ) L\cong \pi^*K$. Choosing a line bundle $K^{1/2}$ we have $L\cong U\pi^*K^{1/2}$ where $\sigma^*U\cong U^*$ and hence  a point in the Prym variety of $S\rightarrow \Sigma$.
The Higgs bundle can be recovered from the spectral data of $S$ together with $U$ by the direct image: $V=\pi_*U$ and $\Phi=\pi_*(U\stackrel{x}\rightarrow U\pi^*K)$.

The map $\det \Phi:{\mathcal M}\rightarrow H^0(\Sigma,K^2)$ defines the integrable system and the generic fibre is the Prym variety. 

\subsection{Symplectic representations of $SL(2,\C)$}
We shall consider the irreducible symplectic representations of $SL(2,\C)$. These are the odd symmetric powers $S^m\C^2$ of the basic representation on $\C^2$ or equivalently the polynomials in one variable of odd degree $m$ under the action 
$$p(z)\mapsto (cz+d)^mp((az+b)/(cz+d)).$$
This is a classical situation and the moment map appears in the literature in different terms. It is a quadratic function with values in the $3$-dimensional Lie algebra $\lie{sl}(2,\C)$ and appears as  a quadratic covariant $C_{2,2}$ in \cite{Ell} or as a net of quadrics in \cite{edge}. In the case of cubic polynomials a modern symplectic treatment is \cite{SS}. It is the invariant homomorphism $S^m\C^2\otimes S^m\C^2\rightarrow S^2\C^2\cong \lie{sl}(2,\C)$ given by multiple contractions with the skew form on $\C^2$.

Here is one more approach to the moment map. Set $m=2k-1$ and write $p(z)=a_0z^m+a_1z^{m-1}+\dots+a_m$, then the invariant symplectic form  on the space $S^m\C^2$   is up to a constant 
$$\omega=\sum_{\ell=0}^{k-1}(-1)^{\ell}\ell!(m-\ell)!da_i\wedge da_{m-i}.$$
  Inside the projective space $\PP(S^m\C^2)$ is the rational normal curve $C$ of perfect powers $(z-\alpha)^m$. The symplectic orthogonal of $p\in S^m\C^2$ is a hyperplane which meets $C$ in $(z-\alpha_i)^m$ where the $\alpha_i$ are the roots of $p(z)=0$. But $p\in p^{\perp}$ so there exist $b_i\in\C$ such that 
$$p(z)=\sum_{1}^mb_i(z-\alpha_i)^m.$$
To evaluate $b_i$ we have $m$ equations 
$$0=p(\alpha_j)=\sum_{i\ne j}b_i(\alpha_j-\alpha_i)^m.$$
For $m$ odd $A_{ij}=(\alpha_j-\alpha_i)^m$ is an odd  skew-symmetric matrix which therefore has a non-trivial kernel, as required. 

The  contraction $S^m\C^2\otimes S^m\C^2\rightarrow S^2\C^2$ then gives 
$$\mu(p)=\sum_{i,j} b_ib_j(\alpha_i-\alpha_j)^{m-1}(z-\alpha_i)(z-\alpha_j)$$
up to a scalar multiple. 

\begin{ex} 

\noindent 1. If $m=1$ then $p(z)=a_0(z-\alpha)$ so $\mu(p)=a_0^2(z-\alpha)^2=p(z)^2$. 

\noindent 2. Take $m=3$ then up to a scalar 
$$p(z)=(\alpha_2-\alpha_3)^3(z-\alpha_1)^3+(\alpha_3-\alpha_1)^3(z-\alpha_2)^3+(\alpha_1-\alpha_2)^3(z-\alpha_3)^3$$
and the formula gives 
$$(\alpha_2-\alpha_3)^3(\alpha_3-\alpha_1)^3(\alpha_1-\alpha_2)^2(z-\alpha_1)(z-\alpha_2)+\dots$$
To rescale and get an expression of degree 2 in $p$ we remove the discriminant factor to get 
\begin{equation}
\mu(p)=a_0^2[\alpha_2-\alpha_3)(\alpha_3-\alpha_1)(z-\alpha_1)(z-\alpha_2)+\dots]
\label{3mom}
\end{equation}
\end{ex}
\vskip .5cm

The discriminant of the quadratic polynomial $\mu(p)$  is  the determinant of the corresponding matrix in $\lie{sl}(2,\C)$ which in our setting defines the quadratic differential and hence the spectral curve.

\subsection{The case $m=1$ }\label{secm1}

The moment map  applied to $u\in \C^2$ is $u\otimes u\in S^2\C^2\cong \lie{sl}(2,\C)$. So given a holomorphic section  $\psi\in H^0(\Sigma, V\otimes K^{1/2})$ we have the nilpotent Higgs field $\Phi=\psi\otimes\psi$. Its kernel is generated by multiples of $\psi$ itself, or otherwise said we have a section of $L^*K^{1/2}$ where $L^*$ is the line bundle defined by $\ker\Phi$ and 
$$0\rightarrow L^*\rightarrow V\rightarrow L\rightarrow 0.$$
The nilpotent cone in the Higgs bundle moduli space consists of the union of  algebraic components each of which is a finite covering of a vector bundle over a symmetric product $S^{2k}\Sigma$ \cite{NH}, together with the moduli space of stable bundles ${\mathcal N}$, where $\Phi=0$. The base of each vector bundle is the divisor of the section of $\Hom(L,L^*K)\cong L^{-2}K$ which is the Higgs field and $k=(g-1)-\deg L$. The fibre consists of the extension classes in $H^1(\Sigma,L^{-2})$ defining $V$. Since $(V,\Phi)$ is stable $\deg L^*<0$, so 
$0< \deg L\le (g-1)$. Then $\deg L^2K> (2g-2)$ so $\dim  H^1(\Sigma,L^{-2})=g-1+2\deg L$. The description just given involves only the line bundle  $L^2$. The $2^{2g}$-fold covering is given by the choice of $L$. When $k>0$ the covering  is connected.

In our case $\psi\in H^0(\Sigma,L^*K^{1/2})$ and $\Phi=\psi\otimes \psi$, so this  translates into the statement that the section of $L^{-2}K$ is $\psi^2$. This means that the divisor is of multiplicity $2$ and so we get the vector bundle restricted to $S^k\Sigma\subset S^{2k}\Sigma$. However now the reduced divisor, together with the choice of square root $K^{1/2}$, which we have fixed from the outset, uniquely determines $L$. Restricted to $S^k\Sigma$ the covering is therefore disconnected and just one component is associated with our choice of $K^{1/2}$, or spin structure.  

The subspace we are considering therefore intersects each component in a space of dimension 
$$k+(g-1)+2((g-1)-k)=3(g-1)-k.$$
When $k=0$, $\deg L=g-1$ and then if $\psi\ne 0$ it is a section of the degree zero line bundle $L^*K^{1/2}$ and so $L=K^{1/2}$.  
The extensions are parametrized by the $3g-3$-dimensional space $H^1(\Sigma, K^*)$. This is precisely half the dimension of ${\mathcal M}$ and so is Lagrangian.

\begin{rmk} The lower dimensional pieces together with the determinant divisor in ${\mathcal N}$ can be shown to form the closure of this $3g-3$-dimensional submanifold, but that is not an issue we shall deal with here, except in the case of genus $2$ and $3$ later on. 
\end{rmk}

\subsection{The case $m=3$ }\label{secm3}
When $m=3$  the moment map is not necessarily nilpotent and we shall consider the generic case where the spectral curve $S$ is smooth. As in Section \ref{spec} The pull-back of  $V$  is expressed as an extension 
$$0\rightarrow L^*\rightarrow \pi^*V\rightarrow L\rightarrow 0.$$
where $L^*$ is  the  eigenspace corresponding to the eigenvalue $x$, a single-valued section of $\pi^*K$ on $S$ such that $\det\Phi=-x^2$. Moreover $L=U\pi^*K^{1/2}$ where $[U]\in \Prym(S)$. 

\begin{prp} \label{m3} Suppose that the spectral curve for $\Phi=\mu(\psi)$, where $\psi$ is a holomorphic section of $S^3V\otimes K^{1/2}$, is smooth. Then $U^3$ is trivial.  
\end{prp}  
\begin{prf} Consider the projective bundle $\PP(V)$ and denote by ${ O}(-1)$ the tautological bundle. Then the section $\psi$ of $S^3V\otimes K^{1/2}$ on $\Sigma$ is the  direct image of a section $s$ of the line bundle $\pi^*K^{1/2}(3)$ on $\PP(V)$. Its divisor is a curve $C$ and the projection $\pi:\PP(V)\rightarrow \Sigma$ represents it as a threefold cover of $\Sigma$. In terms of cubic polynomials the fibre over $a\in \Sigma$ consists of the roots of $\psi(a)$.

Similarly, the Higgs field $\Phi\in H^0(\Sigma, \End_0V\otimes K)= H^0(\Sigma, S^2V\otimes K)$ defines a section $t$ of $\pi^*K(2)$. The line subbundle $L^*\subset \pi^*V$ tautologically defines a section of $\PP(V)$ over $S$ which is the divisor of $t$.  On $S\subset \PP(V)$ the tautological bundle is  $L^*$ and so we have an isomorphism $L^*\cong {O}(-1)$ on $S$.

Now suppose $\Phi=\mu(\psi)$, then from (\ref{3mom}) $\det\Phi$ vanishes when 
$$\mu(\psi)=a_3^2[\alpha_2-\alpha_3)(\alpha_3-\alpha_1)(z-\alpha_1)(z-\alpha_2)+\dots]$$
has a repeated root. Clearly if $\alpha_1=\alpha_2=\alpha$ this is true and the common root is $\alpha$. The quartic invariant $\det\mu(p)$ of a cubic polynomial is its discriminant \cite{edge},\cite{SS} and  in fact it is the only quartic invariant. Thus $C$ is tangential to a fibre of $\pi$ when $\det\Phi=0$ but this is also the ramification of $\pi$ on $S$. This implies that $C$ and $S$ are tangential at these points, or more concretely  $s$ has a double zero on $S$.

Now considering $S$ in the total space of $K$, $\det \Phi$ vanishes when $S$ intersects the zero section $x=0$. Then  the divisor of $s$ restricted to $S\subset \PP(V)$ has a double zero on the divisor of $x$, a section of $\pi^*K$. Hence $\pi^*K^2\cong \pi^*K^{1/2}(3)$  or equivalently $\pi^*K^{3/2}\cong { O}(3)$. Since $L\cong{ O}(1)$ and $L=U\pi^*K^{1/2}$ we deduce that $\pi^*K^{3/2}\cong U^3\pi^*K^{3/2}$ and hence $U^3$ is trivial. 
\end{prf} 

\begin{rmk}  When $U$ itself is trivial then we are in the situation of the canonical section of the integrable system: $V=K^{-1/2}\oplus K^{1/2}$ and $\Phi(u_1,u_2)=(qu_2,0)$ for $q\in H^0(\Sigma, K^2)$. For any simple Lie group these exist \cite{NH1}. They are always Lagrangian, being the fixed point set of an anti-symplectic holomorphic involution. Using this as an identity element, the subvariety described in the proposition consists of the nontrivial 3-torsion points in the abelian variety. This space is connected  as shown by a monodromy calculation in \cite{BS}. Without reference to Higgs bundles we could actually see it is  Lagrangian using local action-angle coordinates by setting the angles to be constant.
\end{rmk}

Proposition \ref{m3} gives a necessary condition but we have not constructed the required section $\psi$ of $S^3V\otimes K^{1/2}$. This we can do for arbitrary odd $m$:
We need to construct now from the isomorphism $L^3\cong \pi^*K^{3/2}$ a section $\psi$ which for which $\mu(\psi)$ is the Higgs field. 
It turns out this works for all odd $m=2k-1$.

Consider the extension  $0\rightarrow L^*\rightarrow \pi^*V\rightarrow L\rightarrow 0$ on $S$. Then there is an induced projection $\pi^*S^mV\otimes K^{1/2}\rightarrow L^m\pi^*K^{1/2}$. 

\begin{prp}   \label{mperiod} Let $S$ be a spectral curve and $U$ a line bundle in the Prym variety.  Suppose $U^{2k-1}$ is trivial and $U^{\ell}$ is non-trivial for all $0<\ell <2k-1$. 
Let $u$ be a  trivialization of $U^{2k-1}$ such that $\sigma^*u=u^{-1}$ and take  the section $ux^k$ of $L^{2k-1}\pi^*K^{1/2}\cong U^{2k-1}\pi^*K^k$ where $x$ is the tautological section of $\pi^*K$. Then this section uniquely extends to the pull-back of a holomorphic section $\psi$ of $S^{2k-1}V\otimes K^{1/2}$ where $V=\pi_*L$. Moreover the resulting Higgs bundle has  $\Phi=\mu(\psi)$ where $\mu$ is the moment map. 
\end{prp}
\begin{rmk} Note that since the moment map is quadratic, $\Phi$ is independent of the choice of $u$.
\end{rmk}
\begin{prf} Defining $V=\pi_*L$ the expression of $\pi^*V$  as an extension 
$$0\rightarrow L^*\rightarrow \pi^*V\rightarrow L\rightarrow 0$$
induces a corresponding decomposition of $\pi^*S^mV$ into a  flag $V_0\subset V_1\subset \dots\subset V_m=\pi^*S^mV$ of subbundles with $V_0=L^{-m}$. Consider the subbundle  $V_k$. In terms of the polynomial $a_0z^m+a_1z^{m-1}+\dots+a_m$ it is the subspace $a_0=a_{1}=\cdots= a_{k-1}=0$, the polynomials which vanish to order $k$ at $\infty$. The subspace $V_k$  is maximally isotropic with respect to the symplectic form on $\pi^*S^mV$.  

It follows that $\pi^*S^mV$ is an extension
\begin{equation}
0\rightarrow V_k\rightarrow \pi^*S^mV\rightarrow V_k^*\rightarrow 0
\label{ext}
\end{equation}

However $\pi^*V$ is also an extension 
$$0\rightarrow \sigma^*L^*\rightarrow \pi^*V\rightarrow \sigma^*L\rightarrow 0$$
and this splits the first extension outside the divisor $D$ given by $x=0$. The extension class is then supported on $x=0$. In Dolbeault terms 
we have a section  of $L\sigma^*L^*$ on $x=0$, extend it holomorphically to a neighbourhood and then everywhere using bump functions to a $C^{\infty}$ section $s$. The extension class in $H^1(S,L^{-2})$ is $[\bar\partial s/x]$. The extension  (\ref{ext}) is for the same reason of the same form.

Dolbeault representatives for an extension $V$ presuppose a $C^{\infty}$-splitting and a $\bar\partial $-operator $\bar\partial +\alpha$ where $\alpha\in \Omega^{01}(\Sigma, L^{-2})$ is the extension form.  So a holomorphic section of $\pi^*S^mV\otimes K^{1/2}$ can be written  as $a_0z^m+a_1z^{m-1}+\dots+a_m$ where 
$$\bar\partial a_0=0,\qquad \bar\partial a_1+m\alpha a_0=0,\qquad \bar\partial a_2+(m-1)\alpha a_1=0, \dots$$
and $a_0\in \Omega^0(\Sigma, L^mK^{1/2}), a_1\in \Omega^0(\Sigma, L^{m-2}K^{1/2}),$ etc.
 
Take $a_0=ux^k$ which is holomorphic. Then 
$$m\alpha a_0=mu\bar\partial s x^{k-1}=\bar\partial (musx^{k-1})$$
so we can take $a_1=-musx^{k-1}$ to get a holomorphic section of $V_1^*\otimes K^{1/2}$. This vanishes to order $(k-1)$ on $D$. Any two extensions of $a_0$ differ by an element of $H^0(S, L^{m-2}\pi^*K^{1/2})$ and if they  are divisible by $x^{k-1}$ then this is  equivalent to sections of $L^{m-2}\pi^*K^{-k+3/2}=U^{m-2}$. If $U^{m-2}$ is non-trivial then such a section is zero and there is a unique extension. 
 Repeating, we get a uniquely determined section $v$ of $V_k^*\otimes \pi^*K^{1/2}$ which vanishes with multiplicity one on $D$.

To extend to $\pi^*S^mV$ we observe that the extension class in $H^1(S, \Hom(V_k^*,V_k))$ for (\ref{ext}) is again supported on $D$ and since $v$ vanishes on $D$ it extends. In this case any two extensions differ by $v'\in H^0(S,V_k)$. But we have 
 $$0\rightarrow V_{\ell}\rightarrow V_{\ell+1}\rightarrow L^{-m+2\ell}\rightarrow 0$$
and  $L^{-m+2\ell}\pi^*K^{1/2}\cong U^{-m+2\ell}\pi^*K^{-k+1+\ell}$ has negative degree for $\ell<k-1$. This means $H^0(S,V_{\ell}\pi^*K^{1/2})\cong H^0(S, V_{\ell +1}\pi^*K^{1/2})$. Since  $V_0\pi^*K^{1/2}=L^{-m}\pi^*K^{1/2}$ has negative degree all sections up to $H^0(S, V_{k-1}\pi^*K^{1/2})$ vanish. 
 Finally  $H^0(S,V_k\pi^*K^{1/2})$ injects into $H^0(S,U^{-1})$ which is zero since $U$ is nontrivial, and so there is a unique extension. 
 
 This gives us a section  $v''$ of $\pi^*S^mV\otimes K^{1/2}$ on $S$. But by choosing a compatible trivialization $u$ we see that $\sigma$ takes $ux^k$ to $u^{-1}x^k$ and by uniqueness $\sigma^*v''=v''$. Since $\pi^*V$ is pulled back from $\Sigma$, $v''$ descends to a section $\psi$ of $S^mV\otimes K^{1/2}.$

It remains to show that the Higgs field is $\mu(\psi)$. Now up to a factor the moment map is the image of $\psi\otimes\psi$ (which is the moment map for $Sp(2k,\C)$) in $S^2V\otimes K$ under the Clebsch-Gordan decomposition of $S^mV\otimes S^mV$. This is a quadratic polynomial $b_0z^2+b_1z+b_2$ where, for some nonzero coefficients $c_i,c_i',c_i''$ 
$$b_0=\sum_{i=0}^{k-1} c_ia_ia_{m-i-1},\qquad b_1=\sum_{i=0}^{k-1} c_i'a_ia_{m-i},\quad b_2=\sum_{i=1}^{k} c''_ia_ia_{m-i+1}.$$
We are not interested in the precise values here but the reader may find  the examples $m=3$ and $5$ in \cite{SS},\cite{edge},\cite{Sal}.

On the divisor $D$ we have $a_0=a_1=\dots =a_{k-1}=0$ and so $b_0=b_1=0$ and $b_2=c''_ka_k^2$. So restricted to the divisor $D$, $\pi^*\mu(\psi)\in H^0(S,\pi^*\End V\otimes K)$ is nilpotent and preserves $L^*\subset V$. This means in particular that $\det\mu(\psi)\in H^0(\Sigma,K^2)$ vanishes on $q=0$ (the spectral curve has equation $x^2-q=0$).

If it is not identically zero then it is a Higgs field with spectral curve $S$. Let $(L')^*$ be the corresponding eigenspace bundle  in $\pi^*V$. Since $\pi^*\mu(\psi)$ preserves $L^*$ on the divisor $D$ we have a holomorphic section of $L'L$ which vanishes on $D$. But $\deg (L'L)=\deg \pi^*K$ so the homomorphism is either zero in which case $L'=L$ or $L'=L^*\pi^*K=\sigma^*L$. Up to a scalar which can be fixed by the definition of $\mu$, this is the required Higgs field $\Phi$.

If $\det\mu(\psi)$ vanishes identically then we have a nilpotent section of $\End V\otimes K$  and a kernel $M\subset V$ on $\Sigma$. Pulling back to $S$ we again have a homomorphism from $\pi^*M$ to $L$ which vanishes on $D$. This cannot be identically zero for then $\pi^*M=L^*$ so that $\Phi$, the Higgs field defined by $L$, has a unique eigenspace. But this means that $\Phi$ is nilpotent and $q=0$. We deduce that  $\pi^*M^*L\pi^*K^*$ has a non-zero section and so 
$$0\le -\deg \pi^*M+\deg L-\deg \pi^*K^*=-2\deg M-2(g-1).$$
However, for a nilpotent Higgs field with kernel $M\subset V$ the Higgs field is a section of $M^2K$ but the inequality above gives  $\deg M^2K\le 0$. It follows that the section must be non-vanishing and $M\cong U'K^{-1/2}$ for a line bundle on $\Sigma$ with $(U')^2$ trivial. But then $\pi^*M^*L\pi^*K^*\cong U(U')^*$  and the non-zero section gives $U'\cong U$ and so $U^2$ is trivial. Since also 
$U^{2k-1}$ is trivial this means that $U$ is trivial which is not true by assumption. 
\end{prf}

\begin{rmk} Although this constraint produces a Lagrangian submanifold it is not obvious to the author that the condition $U^{2k-1}$ trivial is always necessary for $m>3$, The geometry of the quartic invariant, which linked us to the spectral curve for $m=3$, is not so clear in general \cite{edge}. 
\end{rmk}

 \subsection{An exotic example}
 The following case (and its outcome) was suggested to the author by D.Gaiotto.  Take the group $G^c=SL(2,\C)\times SL(2,\C)\times SL(2,\C)$ acting on $U_1\otimes U_2\otimes U_3$, the tensor product of the three 2-dimensional representations.  Since each is symplectic, the tensor product of an odd number is  symplectic. We shall also use the fact that the product of two has an orthogonal structure so that for example  $SL(2,\C)\times SL(2,\C)$ acts on  $U_2\otimes U_3$ as $SO(4,\C)$.
 
 Any Higgs field is of the form 
$$\Phi=\phi_1\otimes 1\otimes 1+1\otimes \phi_2\otimes 1+1\otimes 1\otimes \phi_3$$
and we must obtain those given by $\mu(\psi)$ for some vector $\psi\in U_1\otimes U_2\otimes U_3$. This involves the invariant contraction 
$$S^2(U_1\otimes U_2\otimes U_3)\rightarrow S^2U_1\oplus S^2U_2\oplus S^2U_3.$$

Regard $U_2\otimes U_3$ as having an invariant inner product and then 
$\phi_1$ is obtained by symmetrizing on $U_1$ and taking the inner product on $U_2\otimes U_3$.
Choose a  basis $u_1,u_2$ for $U_1$ and write $\psi=u_1\otimes e_1+u_2\otimes e_2$ then 
$$\phi_1= (e_1,e_1) u_1\otimes u_1+(e_1,e_2)(u_1\otimes u_2+u_2\otimes u_1)+(e_2,e_2)u_2\otimes u_2$$
And then 
\begin{equation}
\tr \phi_1^2= 2\langle u_1,u_2\rangle^2((e_1,e_1)(e_2,e_2)- (e_1,e_2)^2).   
\label{phi1}
\end{equation}
 Now write 
 $$e_1=v_1\otimes w_1+v_2\otimes w_2,\qquad e_2=v_1\otimes (aw_1+bw_2)+v_2\otimes (cw_1+dw_2)$$ in terms of $U_2$ and $U_3$.
 Then
 $$(e_1,e_1)=  2\langle v_1,\,v_2\rangle \langle w_1,w_2\rangle,\, (e_2,e_2)=2(ad-bc)\langle v_1,v_2\rangle \langle w_1,w_2\rangle,  (e_1,e_2)=(a+d)\langle v_1,v_2\rangle \langle w_1,w_2\rangle$$
 and 
\begin{equation}
\tr \phi_1^2=2\langle u_1,u_2\rangle^2\langle v_1,v_2\rangle^2\langle w_1,w_2\rangle^2(4(ad-bc)-(a+d)^2).
\label{form1}
\end{equation}
 Now use the basis $v_1,v_2$ of $U_2$ to write $\psi=v_1\otimes f_1+v_2\otimes f_2$. This gives 
 $$\psi=v_1\otimes (u_1\otimes w_1+u_2\otimes (aw_1+bw_2))+v_2\otimes (u_1\otimes w_2+u_2\otimes (cw_1+dw_2))$$
 and 
 $$(f_1,f_1)=2b\langle u_1,u_2\rangle \langle w_1,w_2\rangle,\, (f_2,f_2)=-2c\langle u_1,u_2\rangle \langle w_1,w_2\rangle,\, (f_1,f_2)=(d-a)\langle u_1,u_2\rangle \langle w_1,w_2\rangle$$
 so, just as in  (\ref{phi1})
 $$\tr\phi_2^2=2\langle v_1,v_2\rangle^2\langle u_1,u_2\rangle^2 \langle w_1,w_2\rangle^2(-4bc -(d-a)^2). $$
 From this and (\ref{form1}) we now obtain  
$\tr\phi_1^2=\tr\phi_2^2$
and similarly 
$$\tr\phi_1^2=\tr\phi_2^2=\tr\phi_3^2.$$
\vskip .5cm
This piece of algebra shows that a Higgs bundle with $\Phi=\mu(\psi)$ for $\psi$ a holomorphic section of $ V_1\otimes V_2\otimes V_3\otimes K^{1/2}$ on $\Sigma$ means that $V_1,V_2,V_3$ are all obtained from line bundles $L_i=U_iK^{1/2}$ on the {\it same} spectral curve $S$.
\vskip .5cm
Consider now the fibre bundle $Q\rightarrow \Sigma$ with fibre $\PP(V_1)\times\PP(V_2)\times \PP(V_3)$. The section $\psi$ defines a section $s$ of the line bundle ${ O}(1,1,1)\otimes \pi^*K^{1/2}$. The line bundles $L_1,L_2,L_3$ on $S$ embed $S$ in $Q$, and on the divisor $D$ given by $x=0$, we have an isomorphism $L_1^*\cong { O}(-1,0,0)$ etc. just as in Section \ref{secm3}. 

The factors $\langle u_1,u_2\rangle^2, \langle v_1,v_2\rangle^2, \langle w_1,w_2\rangle^2$ in $\tr\phi_1^2$ etc. mean that the  hypersurface defined by $s=0$ meets the curve $S$ tangentially at the zeros of $x$, a divisor of $K$, and hence on $S$, $K^{1/2}(1,1,1)\cong K^2$, or equivalently 
$$L_1L_2L_3\cong K^{3/2}.$$

Conversely, we have on $S$ a filtration of $V_1\otimes V_2\otimes V_3$ and so a projection to $L_1L_2L_3$. If $L_1L_2L_3\cong K^{3/2}$, then  $L_1L_2L_3K^{1/2}\cong K^2$ and we can take the section $x^2$. This extends as in Proposition \ref{mperiod}. 

If ${\mathcal M}$ denotes the $SL(2,\C)$-moduli space then the Lagrangian in ${\mathcal M}^3$ maps to the $(3g-3)$-dimensional diagonal in $H^0(\Sigma,K^2)^3$ and the generic fibre is the $2(3g-3)$-dimensional subgroup of $\Prym(S)^3$ given by $[L_1L_2L_3]=1$.
\section{Genus 2}
\subsection{Moduli spaces}
The case of a curve of genus $2$ offers the opportunity  to see some of these Lagrangians more concretely because the moduli space of (S-equivalence classes of) semistable  bundles ${\mathcal N}$ has a concrete description as $\PP^3$ \cite{NR}. Moreover the integrable system as a map from $T^*\PP^3$ has been written down explicitly \cite{GT},\cite{Prev}. 

If $\Sigma$ has genus $2$ it is a hyperelliptic curve $y^2=p(z)$ where $p$ is a polynomial of degree $6$ with distinct roots $z_1,\dots,z_6.$ Invariantly we can regard $p(z)$ as a holomorphic section of ${ O}(6)$ on $\PP^1$ and $y$ as the tautological section of the pull-back of  ${O}(3)$ on the total space, realizing $\Sigma$ as a divisor of $\pi^*O(3)$. Following \cite{NR}  a rank $2$ semistable bundle $V$ with $\Lambda^2V$ trivial is determined by its  one-parameter family of line subbundles $L^*$ of degree $-1$ and these  define  a $2\Theta$-divisor on $\Pic^{1}(\Sigma)$ and hence  a point in the 3-dimensional projective space $\PP(H^0(\Pic^1(\Sigma), 2\Theta))$. If $V=U\oplus U^*$ where $U$ has degree zero, then $[U]\in \Pic^0(\Sigma)$ maps to the Kummer quartic surface in $\PP^3$, where $U$ and $U^*$ are identified since $V$ does not distinguish them. This is the semistable locus. Line bundles $U$ such that $U^2$ is trivial define $16$ distinguished points on it, the singularities of the Kummer surface. 

Unfortunately there is still no uniform description for ${\mathcal M}$, only a stratification depending on the type of the holomorphic structure of $V$ \cite{NH}.
\subsection{Lagrangians for $m=1$}
From Section \ref{secm1} the Lagrangian has an open set $S$ consisting of $V$ expressed as an extension $0\rightarrow K^{-1/2}\rightarrow V\rightarrow K^{1/2}\rightarrow 0$ and such extensions are parametrized by the $3$-dimensional space $H^1(\Sigma, K^{-1})\cong \C^3$. In \ref{secm1} we saw that lower-dimensional pieces are  given by extensions with  $0\le \deg L\le g-1$, but then  $L=K^{1/2}$ is the only one in genus $2$ other than the case $\Phi=0$. Polystable Higgs bundles are of the form $U\oplus U^*$ with $\Phi$ diagonal but then to be in the nilpotent cone we must have $\Phi=0$. So the closure of $S$ consists of points where $\Phi=0$.

The Higgs bundle moduli space itself has a determinant divisor which contains $S$ and therefore also its closure, so the determinant divisor in $\PP^3$ must contain all the limit points.  Conversely if $[V]$ is in this divisor there exists a holomorphic section $\psi$ of $V\otimes K^{1/2}$. Stability for  $V$ implies stability for the Higgs bundle $(V,\mu(\psi))$ and so all points of the determinant divisor of $\PP^3$ are limit points of $S$. 

We know in general that the cohomology class of the determinant divisor of ${\mathcal N}$ generates the cohomology and so in the case of $\PP^3$ this is a projective plane. To describe it we consider $V$ as an extension $0\rightarrow L^*\rightarrow V\rightarrow L\rightarrow 0$ where $L$ has degree $1$. If the family of subbundles $L^*$ contains $K^{-1/2}$ then we have a section $\psi$ of $V\otimes K^{1/2}$ and then considering $\lambda \psi$, $[V]$ is also clearly a limit of points in $S$. The $2\Theta$-divisors on $\Pic^1(\Sigma)$  that pass through our given square root $[K^{1/2}]$  form a plane in $\PP^3$. The $16$ planes formed this way are classically known as {\it tropes} and for convenience we shall retain this name. 

When $V=U\oplus U^*$ is a direct sum of degree zero line bundles, its class lies in the trope if $UK^{1/2}$ or $U^*K^{1/2}$ has a section. By Serre duality this occurs at the same time  and is when $UK^{1/2}$ lies in the theta-divisor of $\Pic^1(\Sigma)$ which is a copy of the curve $\Sigma$ itself since the unique section vanishes at a single point. Thus a trope meets the Kummer surface in the quotient of $\Sigma$ by the involution $\sigma$ which is $\PP^1$ embedded as a conic in $\PP^3$: the trope is tangential to the quartic surface along the conic. It meets a singularity of the Kummer surface where $U^2$ is trivial and $UK^{1/2}$ has a section. Since $(UK^{1/2})^2\cong K$ these are the six odd theta characteristics. In this case the section vanishes at one of the six ramification points $p_1,\dots,p_6$ of $\Sigma\rightarrow \PP^1$. We obtain then the classical configuration of $16$ points and $16$ planes in $\PP^3$, each plane meeting $6$ of the points.

Note that if $U^2$ is non-trivial then  sections $s_1$ of $UK^{1/2}$ and $s_2$ of $U^*K^{1/2}$ vanish at distinct points $x,y$. This means that the section $(s_1,s_2)$ of $(U\oplus U^*)\otimes K^{1/2}$ embeds $K^{-1/2}$ in $U\oplus U^*$ and describes it as an extension as above. On the other hand if  $U^2$ is trivial then $U\oplus U^*$ does not contain $K^{-1/2}$. Instead we should take an {\it S-equivalent} bundle which is a non-trivial extension $0\rightarrow  U\rightarrow V\rightarrow U\rightarrow 0$ and then the section $s_1\in H^0(\Sigma, UK^{1/2})$ lifts to an embedding of $K^{-1/2}$.

We can now identify a compact Lagrangian in the Higgs bundle moduli space ${\mathcal M}$ which is the closure of $S$.  Consider $H^1(\Sigma,K^{-1})\oplus \C\cong \C^4$. For each extension class $[a]\in H^1(\Sigma,K^{-1})$ we have  $0\rightarrow K^{-1/2}\rightarrow V\rightarrow K^{1/2}\rightarrow 0$ and a canonical Higgs field $\phi_{[a]}=1\in H^0(\Sigma, \Hom(K^{1/2},K^{-1/2}\otimes K)\in H^0(\Sigma, \lie{g}\otimes K)$. Then $([a],\lambda)\mapsto (V,\lambda \phi_{[a]})$ defines a map from the projective space  $\PP^3=\PP(H^1(\Sigma,K^{-1})\oplus \C)$ to ${\mathcal M}$. When $\lambda\ne 0$ this is $S$ and when $\lambda=0$ we have the trope.

\begin{rmk}  If we remove the origin in $H^1(\Sigma, K^*)$ which is $V=K^{-1/2}\oplus K^{1/2}$, unstable as a holomorphic bundle, the Lagrangian is the total space of the line bundle ${ O}(1)$ over $\PP^2$. This is not the conormal bundle of the determinant divisor as discussed above, which is ${O}(-1)$. However, the moduli space of {\it stable} bundles is $\PP^3$ with the Kummer quartic surface removed. The ones for which $V\otimes K^{1/2}$ has a section is the complement of a conic in the plane. The section of ${ O}(2)$ defining the conic gives an isomorphism ${ O}(-1)\rightarrow { O}(1)$ outside the conic so the conormal bundle of $\PP^2$ is also isomorphic to ${ O}(1)$ on this open set. 
\end{rmk}

\subsection{The determinant divisor for $m=3$}
We gave in Section \ref{secm3} a geometric description -- the 3-torsion points in the Prym variety -- for the Lagrangian corresponding to the associated symplectic bundle $S^3V$. This was a natural description from the integrable system viewpoint. We want now to see it from the point of view of the conormal bundle in $T^*\PP^3$, that is to describe the determinant divisor as a surface $Y$ in $\PP^3$. The degree of this surface in general can be determined from the ratio of the second Chern class of $S^mV$ in a family to that of $V$ (this is what occurs in the Grothendieck-Riemann-Roch formula). If $m=2k-1$ this ratio is $k(4k^2-1)/3$ and so for $m=3$ it gives $10$. 

Rather than give an explicit polynomial for $Y$ we shall consider the intersection  with the tropes. Since a degree $10$ surface is uniquely determined by its restriction to $16$ planes we lose no information. 

First consider the polystable bundles $V=L\oplus L^*$, $L$ of degree $0$. Then 
$$S^3V=L^3\oplus L\oplus L^{-1}\oplus L^{-3}$$
and so $S^3V\otimes K^{1/2}$ has a section in particular when $LK^{1/2}$ has a section. Thus $Y$ intersects the standard  trope defined by our choice of $K^{1/2}$ in a reducible curve, one component of which is a double conic. The residual part of the curve is of degree $6$. 

Before describing the sextic curve, we make one useful observation about the trope and its conic. Firstly, from our description there is a parametrization of the plane which is independent of the choice of square root $K^{1/2}$, namely as $\PP(H^1(\Sigma, K^{-1}))$. But we can go further. In defining $\Sigma$ as $y^2=p(z)$ in the total space of ${O}(3)$ over $\PP^1$, the canonical bundle $K\cong {O}(1)$ and so $H^1(\Sigma, K^{-1})\cong H^1(\Sigma, {\mathcal O}(-1))$.

From the exact sequence of sheaves ${\mathcal O}(-4)\stackrel{y}\rightarrow {\mathcal O}(-1)$ we have 
$$0\rightarrow H^0(D, {\mathcal O}(-1))\rightarrow H^1(\Sigma, {\mathcal O}(-4))\rightarrow H^1(\Sigma, {\mathcal O}(-1))\rightarrow 0$$
where $D$ is the divisor $y=0$ consisting of the six ramification points $p_1,\dots, p_6$. The 3-dimensional space $H^1(\PP^1,{\mathcal O}(-4))$ injects into $H^1(\Sigma, {\mathcal O}(-4))$ as the $\sigma$-invariant subspace and the connecting homomorphism is anti-invariant since $\sigma(y)=-y$. It follows that each element of the 3-dimensional space $H^1(\Sigma, { O}(-1))$ is uniquely of the form $y\beta$ for some $\beta\in H^1(\PP^1,{\mathcal O}(-4))$. But  under the action of $SL(2,\C)$, $H^1(\PP^1,{\mathcal O}(-4))$ is the standard representation of $SO(3,\C)$ with an invariant inner product. The null vectors define a conic in the projective space.
\begin{lemma}\label{coniclemma}
 This conic  is the same conic as the intersection of the trope with  the Kummer surface.
\end{lemma}
\begin{prf}
Consider $ H^1(\PP^1,{\mathcal O}(-4))$ as  $\lie{sl}(2,\C)$ and a null vector $\beta$ as a nilpotent element. Then acting on $\C^2$ it has a nontrivial kernel. Since $H^0(\PP^1,{\mathcal O}(1))$ is the 2-dimensional representation, there exists a section $s$ of ${O}(1)$ such that $\beta s=0$ in $ H^1(\PP^1,{\mathcal O}(-3))$. Since $K\cong {O}(1)$ and the $g=2$-dimensional space of sections is pulled back from $\PP^1$ this means there is a section $s$ of $K$ such that $y\beta s=0$.

Now look at the  extension defined by $\alpha=y\beta\in H^1(\Sigma,K^{-1})$
$$0\rightarrow 1\rightarrow V\otimes K^{1/2}\rightarrow K\rightarrow 0.$$ The condition  $y\beta  s=0$ means that $s$ defines a section of $K$ which lifts to a section $\psi_1$ of $V\otimes K^{1/2}$. If $\psi_0$ is the section of $V\otimes K^{1/2}$ from the inclusion of the trivial bundle $1$ then $\langle \psi_0,\psi_1\rangle=s$, a section of $K$ which therefore vanishes at points $x$ and $\sigma(x)$. So at $x$ $\psi_1$ is a multiple of the non-vanishing section $\psi_0$  hence $\psi_1-\lambda \psi_0$ vanishes at $x$ and embeds a line bundle $LK^{1/2}\cong{ O}(x)$ in $V\otimes K^{1/2}$. Similarly (if $\sigma(x)\ne x$) another combination $\psi_1-\mu \psi_0$ gives $L^*K^{1/2}$ and $V=L\oplus L^*$.The generic points on the null cone therefore define extensions which correspond to the intersection of the  trope with the Kummer surface, the double conic. 

The six special points $z_i$ correspond as noted above to $V$ which is not a direct sum but an extension, so now we have an extension in two ways 
$$0\rightarrow U\rightarrow V\rightarrow U\rightarrow 0\qquad 0\rightarrow K^{-1/2}\rightarrow V\rightarrow K^{1/2}\rightarrow 0$$
where $U$ coincides with $K^{-1/2}$ at the zero set of a section $s_i$ of $UK^{1/2}$,  the points $p_i\in\Sigma$. Then the class $[\alpha]\in H^1(\Sigma, K^{-1})$ is the image of the isomorphism  $a:U\cong K^{1/2}$ on $p_i$ in the long exact sequence of 
$$0\rightarrow {\mathcal O}(K^{-1})\stackrel{s_i}\rightarrow  {\mathcal O}(K^{-1/2}U)\rightarrow {\mathcal O}_{p_i}(K^{-1/2}U)\rightarrow 0.$$
This means $[\alpha s_i]=0$ and so $[y\beta s_i^2]=0$ and then $[u_i\beta]=0\in H^1(\PP^1,{\mathcal O}(-3))$ where $u_i=s_i^2\in H^0(\PP^1,{\mathcal O}(1))$ vanishes at $z_i$.  Hence $[\beta]$ has a nonzero kernel and again lies on the null conic. 

\end{prf}

 \subsection{The standard trope}

This approach enables us to describe plane curves invariantly associated to the sextic polynomial $p(z)$ defining the hyperelliptic curve. Using the fact that $S^6\C^2$ as a representation of $SL(2,\C)$ is isomorphic to the 7-dimensional space of harmonic cubic polynomials on $\C^3$ we have the following.

\begin{prp} \label{tropeA}Let $\phi(x_1,x_2,x_3)$ be the harmonic polynomial of degree $3$ corresponding to the sextic $p(z)$ defining the genus $2$ curve. Then the intersection of the determinant divisor for $S^3V\otimes K^{1/2}$ with the trope defined by $K^{1/2}$ is the reducible degree $10$ curve with components  the null conic with multiplicity $2$ and the sextic defined by  
$$(x,x)^2\Delta^2\phi^2-16(x,x)\Delta \phi^2-3456 \phi^2=0$$
where $\Delta$ is the Laplacian. 

It  meets  the conic tangentially at the six singular points of the Kummer surface which lie in this trope. \end{prp}

\begin{prf}
We give here the method and leave details of the evaluation until the Appendix. 

For the last part, note that setting $(x,x)=0$ in the equation gives $\phi(x)^2$ and so the sextic meets the conic tangentially at $6$ points given by the intersection of the cubic curve $\phi(x)=0$. These correspond to the six zeros of the sextic $p(z)=0$ under the identification of the null conic with $\PP^1$. 

The bundle $V$ is defined by an extension class in $H^1(\Sigma,K^{-1})$ and, representing this by $\alpha\in \Omega^{01}(\Sigma, K^{-1})$, to find a section of $S^3V\otimes K^{1/2}$ 
we want to solve (as in Proposition \ref{m3} and after some rescaling to simplify notation) 
$$\bar\partial a_0=0,\qquad \bar\partial a_1+\alpha a_0=0,\qquad \bar\partial a_2+\alpha a_1=0, \qquad \bar\partial a_3+\alpha a_2=0$$
where $a_i$ is a $C^{\infty}$ section of $K^{2-i}$. The first equation says that   $a_0$ is holomorphic and given  $a_0\in H^0(\Sigma, K^2)$ the second equation imposes a constraint $[\alpha a_0]=0\in H^1(\Sigma, K)\cong \C$. If this constraint is satisfied then $a_1$ exists and any two choices differ by a holomorphic section $s$ of $K$. The class $[\alpha a_1]\in H^1(\Sigma, {\mathcal O})$ then changes by adding $[\alpha s]$ so if $[\alpha]:H^0(\Sigma,K^2)\rightarrow H^1(\Sigma, {\mathcal O})$ is invertible we can make that choice uniquely. But we have seen in Lemma \ref{coniclemma} that this is the condition that $[\alpha]$ is non-null with respect to the inner product. We make that assumption since we want to determine the sextic curve, the component different from the conic. 

The last  constraint is $[\alpha a_2]=0\in H^1(\Sigma,K^{-1})=H^0(\Sigma, K^2)^*$. The choice in $a_2$ from $\bar\partial a_2+\alpha a_1=0$ is addition of a constant which means the obstruction lies in the cokernel of $[\alpha]:H^0(\Sigma, {\mathcal O})\rightarrow H^1(\Sigma, K^{-1})$. We therefore have a homomorphism mapping $a_0$ in the kernel of $[\alpha]:H^0(\Sigma,K^2)\rightarrow H^1(\Sigma, K)$ to the cokernel of the map $[\alpha]:H^0(\Sigma, {\mathcal O})\rightarrow H^1(\Sigma, K^{-1})$. This is a linear transformation from a $2$-dimensonal vector space to its dual and we want its determinant to vanish. It is homogeneous of degree $3$ in $\alpha$ and so defines a sextic. 

This  linear map is in fact defined by a symmetric form on the kernel of $\alpha$, for suppose we have $(a_0,a_1,a_2)$ and $(a_0',a_1',a_2')$ satisfying these relations, then 
$$\int_{\Sigma} (\alpha a_2) a'_0=\int_{\Sigma}  a_2 (\alpha a_0')=  - \int_{\Sigma}  a_2 (\bar\partial  a_1') =  \int_{\Sigma}  (\bar\partial a_2)   a_1'=- \int_{\Sigma}  \alpha a_1  a_1'$$
using Stokes' theorem, and the right hand side is symmetric. Moreover, it only depends on the components $a_1,a_1'$.

These are statements about the hyperelliptic curve $\Sigma$ but they can be translated into calculations on $\PP^1$. As remarked above,
 $\alpha=y\beta$ where $\beta$ is the   pull-back of a representative for $H^1(\PP^1,{\mathcal O}(-4))$.  Sections of $K$ are pulled back from the two-dimensional space of sections of ${ O}(1)$ on $\PP^1$ and the $3g-3=3$-dimensional space of sections of $K^2$ is pulled back from ${ O}(2)$. 
 
So we can solve $\bar \partial a_1+\alpha a_0=0$ by taking $b_0\in H^0(\PP^1,{\mathcal O}(2))$ such that $[\beta b_0]=0$ and then $-\beta b_0=\bar\partial b_1$ for $b_1$ a $C^{\infty}$ section of ${ O}(-2)$ and setting $a_1=-yb_1$ thought of as a section of ${ O}(1)=K$ on $\Sigma$. The quadratic form in these terms is a multiple of 
$$ \int_{\Sigma}  \alpha a_1^2=  \int_{\Sigma}  y^2\alpha b_1^2= \int_{\Sigma}  y\beta p  b_1^2=\int_{\PP^1}  \beta p  b_1^2$$
since $y:H^1(\PP^1,{\mathcal O}(-2))\rightarrow H^1(\Sigma,K)$ is an isomorphism. In the Appendix we carry out this computation. 
\end{prf}

\subsection{The general trope}
The previous discussion concerns the  particular trope singled out by the choice of the spin structure $K^{1/2}$. Considering the intersection of the determinant divisor for $S^3V\otimes K^{1/2}$ with the trope  for $UK^{1/2}$ where $U^2$ is trivial is equivalent to  asking for a section of $S^3V\otimes UK^{1/2}$ and its intersection with the standard trope. The $15$ non-trivial line bundles $U$ are parametrized by the $15$ choices of pairs of the  ramification points $p_1,\dots,p_6$ of $\Sigma$: the line bundle $UK$ has a unique section $s$ whose divisor is $p_i+p_j$. Making such a choice picks out two zeros $z_i,z_j$ of $p(z)$.  Then $p=q(z)r(z)$ for a quadratic polynomial $q(z)=(z-z_i)(z-z_j)$  and a quartic polynomial $r(z)$. 

\begin{prp} \label{tropeB}The intersection of the  determinant divisor for $S^3V\otimes UK^{1/2}$ with the trope defined by $K^{1/2}$ is a reducible curve with components  a quartic $C_4$ and a sextic $C_6$.
\end{prp}
\begin{prf}
Again we give the method, analogous to the previous case, taking $p_1$ and $p_2$ to define $U$. Given $a_0\in H^0(\Sigma,UK^2)$ we   now have $H^1(\Sigma, KU)=0$ as $U$ is non-trivial, so there exists $a_1$ such that $\bar\partial a_1+\alpha a_0=0$ and $a_1$ is uniquely determined up to the addition of a holomorphic section of $UK$ i.e a multiple of  the section $s$  which vanishes at $p_1,p_2$. The class $[\alpha a_1]\in H^1(\Sigma, U)$ can be killed by adding to $a_1$ a multiple of $s$ so long as $[\alpha s]\ne 0$ in the one-dimensional space  $H^1(\Sigma, U)$, or equivalently if $[\alpha s^2]\ne 0\in H^1(\Sigma,K)$. In that case we find $a_2$ with $\bar\partial a_2+\alpha a_1=0$ and the final obstruction is $[\alpha a_2]\in H^1(\Sigma, UK^{-1})$. As in the previous case this is equivalent to the quadratic form 
$$\int_{\Sigma} \alpha a_1^2v^2$$
for $v\in H^0(\Sigma, UK^2)$ being degenerate. 

The vector space $H^0(\Sigma, UK^2)$ has a natural direct sum decomposition. The section $y$ of ${ O}(3)$ vanishes at the points $p_1,p_2$ and so $y=st$ for a section of $U^*K^{-1}(3)$. But since $K\cong { O}(1)$ and $U^2$ is trivial,  multiples of $t$ give  a  distinguished one-dimensional subspace.  There is a complementary 2-dimensional subspace of the form $su$ where $u\in H^0(\Sigma, K)$.   The equations we need to solve are  linear in $a_0$ and so we can treat the two types separately, giving two distinct components for the determinant divisor. The two intersect if $a_0=0$.

We write $y^2=s^2t^2=qr=p$  and transfer the calculations to $\PP^1$ as follows, beginning with the one-dimensional subspace.

\noindent 1. Take $a_0=t$. Then $\alpha=y\beta=st\beta\in \Omega^{01}(\PP^1,O(-4))$ and 
$$\alpha a_0=yt\beta=sr\beta=\bar\partial (sb_1)$$  
if we can solve $\bar\partial b_1=r\beta$ on $\PP^1$. In that case $a_i=-sb_1$ solves the first equation. But $r\beta\in \Omega^{01}(\PP^1)$  since $r\in H^0(\PP^1,{\mathcal O}(4))$. Given that   $H^1(\PP^1,{\mathcal O})=0$   a solution exists, unique up to adding a constant.

Now
$$\alpha sb_1=y\beta s b_1=qt\beta b_1=\bar\partial (tb_2)$$
if we can solve $\bar\partial b_2=qb_1\beta$ on $\PP^1$. This requires $[qb_1\beta]\in H^1(\PP^1,{\mathcal O}(-2))\cong \C$ to vanish. But if $[q\beta]\in H^1(\PP^1,{\mathcal O}(-2))$ is nonzero, this can be achieved by the choice of constant in $b_1$.  (In the following we shall regard $[q\beta]$ and $[qb_1\beta]$ as complex numbers).

The final condition is the vanishing of the class $[\alpha tb_2]=[srb_2\beta]\in H^1(\Sigma, UK^{-1})$ which is equivalent to 
$$0= [rb_2\beta]\in H^1(\PP^1,{\mathcal O}(-2)).$$
Since $H^1(\PP^1,{\mathcal O}(-2))=H^1(\PP^1,K)$ this is an integral and we can rewrite it as 
$$\int_{\PP^1}rb_2\beta=\int_{\PP^1}b_2\bar\partial b_1=-\int_{\PP^1}\bar\partial b_2 b_1=-\int_{\PP^1}b_1^2 q \beta.$$
This makes sense for any choice of $b$ solving $\bar\partial b=r\beta$ which is linear in $\beta$ and is a homogeneous  expression $C([\beta])$ cubic in $\beta$. But in order to define $b_2$ we  modify this via 
$$b_1=b-\frac{[qb\beta]}{[q\beta]}$$
and this gives 
$$\int_{\PP^1}b^2 q \beta-2\frac{[qb \beta]}{[q\beta]}\int_{\PP^1} qb\beta+\frac{[qb \beta]^2}{[q\beta]^2}\int_{\PP^1} q\beta=C([\beta])-\frac{[qb \beta]^2}{[q\beta]}.$$
The constraint is therefore given by a quartic curve with the following equation
\begin{equation}
{[q\beta]}C([\beta])-{[qb \beta]^2}=0
\label{C4}
\end{equation}
This is  therefore a component of the determinant divisor.

\noindent 2. Now, for the 2-dimensional complementary subspace, consider $a_0$ of the form $su$ for $u\in H^0(\PP^1,{\mathcal O}(1))$. We have 
$$\alpha a_0= s^2tu\beta=tqu\beta=\bar\partial (tb_{1})$$
 But this lies in $H^1(\PP^1,{\mathcal O}(-1))=0$ and moreover since $H^0(\PP^1,{\mathcal O}(-1))=0$ there is a unique $b_{1}$.
Next
$$\alpha t b_{1}=srb_{1}\beta=\bar\partial (sb_{2})$$
can be solved if $[rb_{1}\beta]=0\in H^1(\PP^1,{\mathcal O}(-1))$ which is again true. The final constraint is the class 
$[\alpha sb_{2}]=[tq b_{2}\beta]$ or equivalently 
$$[q b_{2}\beta]\in H^1(\PP^1,{\mathcal O}(-3)).$$
This linear map from $u\in H^0(\PP^1,{\mathcal O}(1))$ to $H^1(\PP^1,{\mathcal O}(-3))$ is represented by a $2\times 2$ matrix which is cubic in $\beta$ and so the vanishing of the  determinant defines a curve of degree $6$. 
\end{prf}

Instead of giving formulas to describe these curves we  derive the intersection properties of $C_4$, $C_6$, the conic and the line joining the points $z_1,z_2$ on the conic. 

\begin{prp} Let $C_2$ be the conic, $C_1$  the line of intersection of the two tropes and $C_4$, $C_6$ the quartic and sextic, components of the determinant divisor. Then
\begin{itemize}
\item
$C_1$ intersects $C_2$ in  $z_1,z_2$, two  of the six distinguished points 
\item
$C_1$ is a bitangent to $C_4$
\item
the two points $C_1\cap C_4$ harmonically separate $z_1,z_2\in C_1$
\item
$C_4$ meets $C_2$ tangentially at the  four points   $z_3,z_4, z_5,z_6.$
\item
$C_6$ meets $C_2$ tangentially at all six points $z_i$.
\end{itemize}
\end{prp}

\begin{prf}

\noindent 1. From Lemma \ref{coniclemma} each of the six points in the trope for $K^{1/2}$ is represented by an extension of the form $s_i^{-1}\bar\partial \varphi$. It lies in the trope for $UK^{1/2}$ if there is a homomorphism from $UK^{-1/2}$ to $K^{1/2}$ which lifts to $V$: equivalently a section $s$ of $UK$ for which $[\alpha s]=0$. But $s$ is unique and vanishes at $p_1,p_2$. Hence if $i=1,2$ $s_i$ divides $s$ and $ss_i^{-1}\bar\partial \varphi=\bar\partial(ss_i^{-1}\varphi)$ so the class is trivial. This means the two points $z_1,z_2$ on the conic also lie in the trope for $K^{1/2}U$ i.e. the line of intersection $C_1$ intersects $C_2$ in $z_1,z_2$. This is the line $[q\beta]=0$.

\noindent 2. In the equation (\ref{C4}) for the quartic substituting $[q\beta]=0$ gives $[qb\beta]^2=0$ and $[qb\beta]$ is quadratic in $\beta$ so $C_4$ meets $C_1$ in two points with multiplicity $2$. The line is therefore a bitangent. 

  We have  on $\PP^1$ the  exact sequence of sheaves 
$$0\rightarrow {\mathcal O}(-4)\stackrel{q}\rightarrow  {\mathcal O}(-2)\rightarrow {\mathcal O}_D(-2)\rightarrow 0$$
where $D$ is the divisor of $q$. From the long exact cohomology sequence this means that the extension class  comes from $H^0(D,{\mathcal O}_D(-2))$: a pair $(a_1,a_2)\in \C^2$.
In Dolbeault terms we  represent $[\beta]$ as before by a  form $\beta=q^{-1}\bar\partial \gamma$ where $\gamma$ is supported in a neighbourhood of $D$ and $\bar\partial\gamma=0$ in a smaller neighbourhood and takes the values $
c_1,c_2$ at the two points $z_1,z_2$ the  zeros of $q$.  Then
$$[qb\beta]=\int_{\PP^1}b\bar\partial \gamma=-\int_{\PP^1}\bar\partial b \gamma=-\frac{1}{2}\int_{\PP^1}rq^{-1}\bar\partial (\gamma^2).$$
and by Stokes' theorem taking $q=(z-z_1)(z-z_2)$ this localizes to 
$$\frac{r(z_1)c_1^2}{2(z_1-z_2)}+\frac{r(z_2)c_2^2}{2(z_2-z_1)}.$$
Now $(c_1,c_2)$ are homogeneous coordinates on the line $[\beta]\in H^1(\PP^1,{\mathcal O}(-4))$ defined by $[q\beta]=0$. Moreover if $c_1=0$, then the section of ${ O}(1)$ which vanishes at $z_1$ is annihilated by the extension class. This  is then a point on the null cone: hence the line $C_1$ intersects the null cone at $z_1,z_2$ with coordinates $(1,0),(0,1)$. The quartic intersects the line then in the two points $\pm \sqrt{r(z_1)/r(z_2)}$ making them harmonically separate $z_1,z_2$.

\noindent 3.  To find the intersection of $C_4$ with the conic $C_2$ we first consider the generic point  where $V=L\oplus L^*$ and $LK^{1/2}$ has a section $s_1$ with zero $x$ and $L^*K^{1/2}$ a section $s_2$ with zero $\sigma(x)$. Assume that $L^2$ is nontrivial so that $x\ne \sigma(x)$.

The bundle $V$  has a projection 
$(u,v)\mapsto us_2+vs_1$ to $K^{1/2}$ which represents it as an extension.  The induced homomorphism from $S^3V= L^3\oplus L\oplus L^{-1}\oplus L^{-3}$ to $K^{3/2}$ is defined by
$(s_2^3,s_2^2s_1, s_2 s_1^2, s_1^3)$. 
The bundle $S^3V\otimes UK^{1/2}$ has a section if $ULK^{1/2}$ or $UL^3K^{1/2}$ has a section. Both cannot hold simultaneously since if so their divisors $x,y$ satisfy $3x-y\sim K$ but if $x\ne y$ there exists a differential on $\Sigma$ with a single simple pole. But the sum of the residues of a differential must be zero so $x=y$, which means $L^2$ is trivial. 

Let $v$ be a section of $L^3K^{1/2}U$ then $a_0\in H^0(\Sigma,UK^2)$ is defined by $s_2^3v$ which has a divisor of the form $3x+y$. Since $t$ has distinct zeros this does not lie on $C_4$.  Similarly if $w$ is a section of $ULK^{1/2}$ then $a_0=ws_2^2s_1$ has a divisor with some  point of multiplicity $2$.

It follows that all intersections with the conic occur among the six distinguished points. But two of these are $z_1,z_2$ which lie on the line $C_1$ and we have already seen that $C_4$ intersects the line at two  different points. Thus $C_4$ must intersect at $z_3,z_4,z_5,z_6$ and by (Galois) symmetry with the same multiplicity 2 at each point.

 At these points we have $a_0=0$. To see this, assume not and  use the representative $u_i^{-1}\bar\partial \gamma$ of $\beta$. Then since $u_i$ for $i>2$ divides $r$ we can take $b_1=u_i^{-1}r\gamma$. Now $[q\beta]\ne 0$ since $z_i$ does not lie on the line joining $z_1$ and $z_2$ so we need the following integral to vanish:
$$\int_{\PP^1}b_1^2q\beta=\int_{\PP^1}u_i^{-2}r^2 \gamma^2 qu_i^{-1}\bar\partial\gamma=\int_{\PP^1} u_i^{-3}qr^2\bar\partial(\gamma^3/3)$$
and this localizes to $z_i$ giving a non-zero multiple of $c^3$, the value of $\gamma$ at that point.

\noindent 4. Finally consider $C_6$. If $[V]=[L\oplus L^*]$  lies on $C_6$ then, with $v$ a section of $UL^3K^{1/2}$,  $s_2^3v=su$ for a section $u$ of $K$ and so has a divisor of the form $p_1+p_2+z+\sigma(z)$. But then $x=p_1=z=\sigma(z)$ or similarly $x=p_2$. This however means that $L=U$ and again $L^2$ is trivial. If  $w$ is a section of $ULK^{1/2}$ then $a_0=ws_2^2s_1$ which has divisor $$2x+\sigma(x)+y=p_1+p_2+z+\sigma(z).$$
Here we could have  $x=z=\sigma(z)$ which means that $x=p_i$ but again this gives $L^2$ trivial. 

In the previous case we saw that  $C_4$  passes through the points $z_i$ for $i>2$ so consider $z_1$ and $a_0=0$ which means that $C_6$ also passes through those points. It remains to consider $z_1,z_2$. Take  $a_0=su_1$.
Using $\beta= u_1^{-1}\bar\partial\gamma$ we see that 
$\bar\partial b_1=qu_1\beta$ is solved by $b_1=q\gamma$. Then $rb_1\beta=rqu_1^{-1}\bar\partial \gamma^2/2$ and  we can take $b_2=\gamma^2 rq u_1^{-1}/2$ since $u_1$ divides $q$. Then 
$$\int_{\PP^1}qb_2\beta u_j=\int_{\PP^1}q^2u_1^{-2}u\bar\partial(\gamma^3/6)=0$$
for all $u\in  H^0(\PP^1,{\mathcal O}(1))$ since $qu_1^{-1}$ is holomorphic.  The quadratic form on $H^0(\PP^1,{\mathcal O}(1))$ therefore vanishes at $z_1$ and $z_2$ and so $C_6$ meets $C_2$ with multiplicity $2$ at  $z_1,z_2$.
\end{prf}

\begin{rmk} The intersections just considered are actually intersections of the 3-dimensional Lagrangians for $m=1$ and $m=3$, which can be seen by considering the filtration $V_0\subset V_1\subset V_2\subset V_3$ of $S^3V$ corresponding to the extension class $[\alpha]\in H^1(\Sigma, K^{-1})$ defining $V$ in the standard trope. Then for $\mu(\psi)$ to be nilpotent with kernel $K^{-1/2}\subset V$ the spinor $\psi$ must lie in $V_1\otimes K^{1/2}\cong V\otimes K^{-1/2}$. This is an extension
$$0\rightarrow K^{-1}\rightarrow V_1\otimes K^{1/2}\rightarrow 1\rightarrow 0.$$
To obtain $\psi$ we need to lift a constant section of the trivial bundle to $V_1\otimes K^{1/2}$ but the obstruction is precisely $[\alpha]$. When this vanishes we get the direct sum $V=K^{-1/2}\oplus K^{1/2}$ and the canonical Higgs bundle. This indeed does have the right section, but it is the intersection of the $m=1$ Lagrangian $\PP^3$ with the canonical section of the fibration. That is, on a general fibre, the zero among the  3-torsion points. The Lagrangian given by the non-trivial elements of order $3$ therefore has a closure which must intersect at $\Phi=0$ on the determinant divisor.
\end{rmk}

 As well as considering the intersection of the degree 10 surface with the tropes we could also have considered its intersection with the Kummer surface, but this is better handled by other methods. It means considering $V=L\oplus L^*$, $\deg L=0$  where $LK^{1/2}$ or $L^3K^{1/2}$  has a section. The map $f:\Pic^1(\Sigma)\rightarrow \Pic^1(\Sigma)$ defined by $f(LK^{1/2})=L^3K^{1/2}$ pulls back the $\Theta$-line bundle to the $9\Theta$-bundle. Since the  map to $\PP^3$ is defined by $2\Theta$ the determinant divisor is in the linear system $20\Theta$ but has a $2\Theta$ component from the $LK^{1/2}$ condition. 
 
 \section{Genus 3}
 When $\Sigma$ has genus $3$, the line bundles of degree $-1$ in $V$ again define a $2\Theta$-divisor in $\Pic^1(\Sigma)$ but in this case the moduli space ${\mathcal N}$ embeds as a quartic hypersurface in $\PP^7$ \cite{NR1} known classically as the Coble quartic. It has a singularity along the semi-stable locus which is the 3-dimensional Kummer variety $\Pic^0(\Sigma)/\Z_2$.  We shall outline here the structure of the Lagrangian for $m=1$.
 
 From Section \ref{secm1} when the section $\psi$ is non-zero, the Lagrangian meets two components of the nilpotent cone. It contains the open set $H^1(\Sigma, K^{-1})\cong \C^{3g-3}=\C^6$ and the total space of a rank $4$ vector bundle $E$ over $\Sigma$. This is where the image  of  the homomorphism $K^{-1/2}\rightarrow V$ defined by $\psi$ lies in a line bundle $L^*$ of degree $-1$. Then $K^{1/2}L^*$ is of degree $1$ and the section vanishes at a unique point $x$. The fibre of $E$ at $x$ is the $3(g-1)-2=4$-dimensional space of extensions $H^1(\Sigma, L^{-2})$.  More concretely, we are looking at $H^1(\Sigma, K^{-1}(2x))$ for $x\in \Sigma$. The dual bundle has fibre $H^0(\Sigma, K^2(-2x))$ which consists of sections whose $1$-jet vanishes at $x$ hence
 \begin{equation}
 0\rightarrow E^*\rightarrow \Sigma\times H^0(\Sigma, K^2)\stackrel{j^1}\rightarrow J^1(K^2)\rightarrow 0.
 \label{Edef}
 \end{equation}
 By using the theory of stable pairs (see e.g. \cite{Brad}) one may see that for $\psi\ne 0$ these two strata form a  smooth manifold and we need to understand the closure of this space which, as discussed above, means adding  the determinant divisor in ${\mathcal N}$. 
 
The  determinant divisor (or Brill-Noether loci in general) forms  the subject of \cite{OPP}. A  starting point  is the fact that for arbitrary $g$ and $L\in \Pic^1(\Sigma)$, the extension
 $$0\rightarrow L^*\rightarrow V\rightarrow L\rightarrow 0$$
 defines an injection of the projective space $\PP(H^1(\Sigma, L^{-2}))$ to the moduli space ${\mathcal N}$ of semistable bundles. Restricting to $\Sigma\subset \Pic^1(\Sigma)$ this is a map of the projective bundle $\PP(E)$ defined above to ${\mathcal N}$ and this  forms part of an explicit description in \cite{OPP} of the determinant divisor ${\mathcal W}$ for $g=3$. One takes the projective space $\PP^5=\PP(H^1(\Sigma, K^{-1}))$  and the bicanonical embedding of the curve $\Sigma\subset\PP^5$. Then the blow-up $\hat\PP^5$ along $\Sigma$ has a well-defined map to ${\mathcal W}$ where for each $x\in \Sigma$ the projectivized normal  space $\PP(N_x)\subset \hat \PP^5$ maps to the embedded $\PP(H^1(\Sigma, K^{-1}(2x))\subset {\mathcal N}$.
 
 This can be viewed in the context of the Morse theory of the function $f=\Vert \Phi\Vert^2$ on ${\mathcal M}$ restricted to our Lagrangian.  Since this is  invariant by the action $\psi\mapsto e^{i\theta}\psi$, $f$ restricts to a moment map on this submanifold. In our case the action has fixed point set at two levels: the maximum of $f$ which is the trivial extension $V=K^{-1/2}\oplus K^{1/2}$ with the canonical Higgs field, a single point in ${\mathcal M}$,  and a lower value $f=c$  which is  where $V=L^*\oplus L$ and is  the zero section of the vector bundle $E$. The downward Morse flow from $f=c-\epsilon$  has limit on ${\mathcal W}$  and so we have a map from the K\"ahler  quotient at this level to ${\mathcal W}$. At the level $f=c+\epsilon$ the quotient is the projective space $\PP(H^1(\Sigma, K^{-1}))$, and following \cite{Th} the birational transformation relating the two is the blow-up of $\Sigma$.

If we  now  take $\PP(H^1(\Sigma, K^{-1})\oplus \C)$ and blow up the submanifold $\Sigma\subset \PP(H^1(\Sigma, K^{-1}))$ in $ \PP(H^1(\Sigma, K^{-1})\oplus \C)$ then we have the projective variety  $\hat \PP^6$ which maps to the closure of the Lagrangian.  

\begin{rmk} 1. Note that the tangent bundle of $\PP^6$ restricted to $\PP^5$  splits as $T\oplus O(1)$ so the projectivized normal bundle of $\Sigma \subset \PP^6$ is $\PP(N\oplus O(1))$ so removing $\PP(N)$ gives the total space of the vector bundle $N(-1)$. But dualizing (\ref{Edef}) we see that this is precisely $E$. Moreover   $\PP(H^1(\Sigma, K^{-1})\oplus \C)\backslash \PP(H^1(\Sigma, K^{-1}))=H^1(\Sigma, K^{-1})$ which gives us the other stratum. Hence $\hat \PP^6\backslash \hat \PP^5$ is visibly the open part of the Lagrangian. 

\noindent 2. 
The image is singular and in particular at the equivalence class of the  unique stable bundle $V$ for which $V\otimes K^{1/2}$ has a 3-dimensional space of holomorphic sections. This arises from 
$$  0\rightarrow V\otimes K^{-1/2} \rightarrow \Sigma\times H^0(\Sigma, K)\stackrel{\mathrm{ev}}\rightarrow K\rightarrow 0$$
\end{rmk}

\section{Appendix}
We give here the method to find the explicit formula for the sextic curve in the standard trope. In principle it can be used to give formulas for $C_4,C_6$ but these are unlikely to be pleasant.
We shall  use convenient Dolbeault representatives for classes in $H^1(\PP^1,{\mathcal O}(-m))$. 

Consider the  form 
$$\frac{\bar z^m}{(1+z\bar z)^n}dz^kd\bar z$$
on $\C$. It behaves  as $z\rightarrow \infty$ like $z^{-n}\bar z^{m-n}z^{2k}\bar z^{2}$ (using $d(1/z)=-dz/z^2)$ so if $n\ge 2k$ and $m\le n-2$ it represents a cohomology class in $H^1(\PP^1,{\mathcal O}(-2k))$ since the canonical bundle of $\PP^1$ is ${ O}(-2)$. For example $m=0,k=1,n=2$ is the volume form which generates $H^1(\PP^1,{\mathcal O}(-2))=H^1(\PP^1,K)$.

To determine whether it is trivial we need to find the obstruction to solving
$$\bar\partial f=\frac{\bar z^m}{(1+z\bar z)^n}d\bar z$$
and interpreting $fdz^k$ as a global section on $\PP^1$. 
Put $y=(1+z\bar z)$ as a function of $\bar z$, regarding $z$ as a constant  and integrate 
$$\frac{(y-1)^m}{z^{m+1} y^n}dy$$
to get
$$\frac{1}{z^{m+1}}\sum_{k=0}^m\frac{(-1)^m}{k-n+1}{m\choose k}y^{k-n+1}.$$

The leading term in the integral is $z^{-(m+1)}y^{m-n+1}dz^k\sim z^{-(m+1)}z^{m-n+1}\bar z^{m-n+1} z^{2k}$ which tends to zero under the previous conditions. However the integral has a pole of order $m+1$  at $z=0$. We can add a term $cdz^k/z^{m+1}$, which is annihilated by $\bar\partial$, to make this regular but this is smooth at infinity only if $m\ge 2k-1$. So
$     0\le m\le 2k-2 $ gives a $(2k-1)$-dimensional space of cohomology classes in $H^1(\PP^1,{\mathcal O}(-2k))$ and hence all. 

\begin{ex} The volume form $$\frac{1}{(1+z\bar z)^2}dzd\bar z$$ integrates to $-1/zydz$. To make this regular we add $dz/z$ but if $\tilde z=1/z$ this is $-d\tilde z/\tilde z$ and has a pole at infinity. So the principal part of the naive integral represents the cohomology class.
\end{ex} 

The curves $C_i$ defined in the previous section  lie in the  projective space of  the 3-dimensional space $H^1(\PP^1,{\mathcal O}(-4))$. A Dolbeault representative for a class $[\beta]\in H^1(\PP^1,{\mathcal O}(-4))$ is the case $k=2$ above. So $m=0,1,2$ and $n\ge 4$ gives us  representatives. Take $n=4$ and then 
$$[\beta]=\frac{v_0+v_1\bar z+v_2\bar z^2}{(1+z\bar z)^4}dz^2d\bar z.$$
The naive integral of this is $$\int \left(\frac{v_0}{zy^4}+\frac{v_1(y-1)}{z^2y^4}+\frac{v_2(y-1)^2}{z^3y^4}\right)dydz^2$$ which gives
$$f=-\frac{1}{z^3} \left(\frac{1}{3y^3}(v_0z^2-v_1z+v_2)+\frac{1}{2y^2}(v_1z-2v_2)+\frac{1}{y}v_2\right)dz^2.$$

The class $[\beta]$ is in the null cone if $[\beta u]=0\in H^1(\PP^1,{\mathcal O}(-3))$ for some  $u\in  H^0(\PP^1,{\mathcal O}(1))$. Writing $u=(w_0+w_1z)dz^{-1/2}$ we have $\bar\partial uf=\beta u$ and so we need the condition that we can subtract off the polar part of  $uf$ as a smooth section of $O(-3)$. Since $dz^{3/2}/z^3$ is regular at infinity we need the vanishing of the $z^{-2}$ and $z^{-1}$ terms in $fu$ which gives 
$$2v_0w_0+v_1w_1=0,\qquad v_1w_0+2v_2w_1=0$$
and the condition for a non-zero $u$ is $v_1^2-4v_0v_2=0$. This is the equation of the null conic. Parametrize it by $v_0=1, v_1=-2t, v_2=t^2$ then the kernel of the action on $H^0(\PP^1,{\mathcal O}(1))$  is $w_0/w_1=t$. Thus the six points are $t=z_1,\dots, z_6$. 
This is illustrative of the method, which we now apply to the standard trope.

The first step in the proof of Proposition \ref{tropeA} is to choose $a_0\in H^0(\PP^1,{\mathcal O}(2))$ such that $[\beta a_0]=0$, where $[\beta]$ is non-null, since we know the null conic lies in the determinant divisor. By the action of $SL(2,\C)$, and knowing the equation of the conic as above, we can take $v_0=v_2=0$, and  $v_1=1$ so that 
$$f=\frac{1}{z^2} \left(\frac{1}{3y^3}-\frac{1}{2y^2}\right)dz^2.$$Then  $[\beta a_0]=0\in H^1(\PP^1,{\mathcal O}(-2))$ means that $a_0=(u_0+u_2z^2)dz^{-1}$. 

To follow the procedure in the proof  we put $u(z)=u_0+u_2z^2$ and the sextic $p(z)=c_0+c_1z+\dots+c_6z^6$ defining the hyperelliptic curve $\Sigma$ and set  
$b_1=(u_0+u_2z^2)f-p_1$
where $p_1$ is the polar part, in this case $p_1=-u_0/6z^2$. So
$$b_1=(u_0+u_2z^2)\frac{1}{z^2} \left(\frac{1}{3y^3}-\frac{1}{2y^2}\right)+\frac{u_0}{6 z^2}=A_0+\frac{1}{z^2}A_2$$

 The criterion then is the degeneracy of the quadratic form in $u_0,u_2$ defined by 
$$\int_{\PP^1}b_1^2p(z)\beta=\int_{\PP^1}(A_0+\frac{1}{z^2}A_2)^2p(z)\frac{ \bar z}{y^4}dzd\bar z=\int_{\PP^1}(A_0+\frac{1}{z^2}A_2)^2p(z)\frac{ y-1}{z y^4}dzd\bar z.$$
Since $y$ is rotation-invariant it is clear that only the coefficients of $z, z^3, z^5$ in $p(z)$ contribute. Evaluating the integral   shows that the quadratic form is singular if 
$$c_3^2=4c_1c_5.$$

Here we fixed $\beta$ and have arrived at a condition on $p$, but applying the action of $SL(2,\C)$ we can revert to an equation in $(v_0,v_1,v_2)$, the constraint on $\beta$. To obtain the formula in Proposition \ref{tropeA} we relate $\phi$ and $p(z)$ by (this is the twistor approach) setting $\phi$ to be the coefficient of $z^6$ in 
$$((x_1+ix_2)+2x_3z-(x_1-ix_2)z^2)^3p(z).$$
Evaluating $$(x,x)^2\Delta^2\phi^2-16(x,x)\Delta \phi^2-3456 \phi^2=0$$ at $(x_1,x_2,x_3)=(0,0,1)$ gives a multiple of $c_3^2-4c_1c_5.$

\vskip 1cm
 {Mathematical Institute,
Radcliffe Observatory Quarter,
Woodstock Road,
Oxford, OX2 6GG}
\end{document}